\documentclass[a4paper,twoside,openright,10pt,twocolumn]{article}
\usepackage{amsmath}
\usepackage{amsfonts}
\usepackage{stackrel}
\usepackage{graphicx}
\usepackage{caption}
\usepackage{subcaption}
\usepackage{enumerate}
\usepackage{fancyhdr}

\newcommand{\note}[1]{\vspace{5 mm}\par \noindent
  \marginpar{\textsc{Note}} \framebox{\begin{minipage}[c]{0.9
        \textwidth} \flushleft \tt #1 \end{minipage}}\vspace{5
    mm}\par}
\renewcommand{\note}[1]{}

\newcommand {\R}{\mathbb R}

\newcommand {\Z}{\mathbb Z}
\newcommand {\N}{\mathbb N}
\newcommand {\e}{\mathrm e}
\renewcommand {\i}{\mathrm i}
\newcommand {\Fix}{\mathrm{Fix}}
\renewcommand {\S}{\mathbf{S}}
\renewcommand {\d}{\mathrm{d}}

\newcommand {\kC}{\mathcal{C}}
\newcommand {\kX}{\mathcal{X}}
\newcommand {\kF}{\mathcal{F}}
\newcommand {\Mat}{\text{Mat}}
\newcommand {\w}{\text{w}}
\newcommand {\re}{\text{Re}}
\usepackage{calc}

\fancypagestyle{myplain}{%
\vspace{-1.5cm}
\fancyheadoffset[LE,RO]{\marginparsep + \marginparwidth}
  \fancyhf{}%
  \chead{--- DRAFT ---}
  \fancyhead[R]{Proceedings of the ASME 2015 International Design Engineering\\Technical Conferences \&\\
Computers and Information in Engineering Conference\\
IDETC/CIE 2015\\
August 2-5, 2015, Boston, Massachusetts, USA\\[1cm]
\large{DETC2015-46203}}%
\fancyfoot[r]{Copyright \copyright~2015 by ASME}}
\fancypagestyle{otherplain}{%
\fancyheadoffset[LE,RO]{\marginparsep + \marginparwidth}
\fancyhf{}%
\fancyfoot[r]{Copyright \copyright~2015 by ASME}}
\fancypagestyle{myplain}{}
\fancypagestyle{otherplain}{}
\title{\vspace{1.5cm}Stability of Hopf bifurcations in time-delayed fully-connected PLL networks}
\author{Diego Paolo Ferruzzo Correa\\dferruzzo@usp.br\\\'Atila Madureira Bueno\\atila@sorocaba.unesp.br\\Jos\'e Roberto Castilho Piqueira\\piqueira@lac.usp.br}
\date{\today}
\begin{document}
\maketitle
\thispagestyle{myplain}
\abstract{%
Dynamics in delayed differential equations (DDEs) is a well studied problem mainly because DDEs arise in models in many areas of science including biology, physiology, population dynamics and engineering. The change of the nature in the solutions in the parameter space for a network of Phase-Locked Loop oscillators was studied in \textit{Symmetric bifurcation analysis of synchronous states of time-delayed coupled Phase-Locked Loop oscillators}. Communications in Nonlinear Science and Numerical Simulation, Elsevier BV, 2014, (on-line version), where the existence of Hopf bifurcations for both cases, symmetry-preserving and symmetry-breaking synchronization was well stablished. In this work we continue the analysis exploring the stability of periodic solutions emerging near Hopf bifurcations in the Fixed-point subspace, based on the reduction of the infinite-dimensional space onto a two-dimensional center manifold. Numerical simulations are presented  in order to confirm our analitycal results. Although we explore network dynamics of second-order oscillators, results are extendable to higher order nodes.\\[0.25cm]

}
\section{Introduction}
\label{sec:intro}
We consider the Full Phase model introduced in~\cite{FerruzzoCorrea2014a} to analyse stability of periodic orbits near Hopf bifurcations emerging in the parameter space $(\mu,\tau)$  at the Fixed-point subspace, which are non degenerative for the case $K>1$. It has been shown that these bifurcations  can cross the imaginary axis in both directions, from the left to the right and from the right to the left. The main approach used for the analysis is the decomposition of the infinite-dimentional space into a 2-dimensional center space corresponding to the imaginary critical simple eigenvalue $\lambda=\pm \i\omega$, $\omega>0$, and an infinite-dimentional space ``orthogonal'' to the first one (the orthogonality condition will be defined below). We will follow closely the theory and procedures presented in~\cite{Kalmar-Nagy2001, Gilsinn2009, Zhao2009, Campbell2009, Guckenheimer1983}.

\section{The Full-phase model}
\label{sec:fullphase_model}
In~\cite{FerruzzoCorrea2014a} the Full-phase model was used to find Hopf bifurcations in the parameter space $(\mu,\tau)$, the general model for a $N$-node, fully-connected, second-order oscillator network is:
\begin{align}
\begin{array}{l}
\ddot\phi_i(t)+\mu\dot\phi_i(t)-\mu-\dfrac{K\mu}{N-1}\displaystyle\sum_{\stackrel[j=1]{j\neq i}{}}^Nf(\phi_i,\phi_j)=0,
\end{array}
\label{eq:fullphasemodel}
\end{align}
$i=1,\ldots,N,$ where:
\begin{align}
\begin{array}{l}
f(\phi_i,\phi_j)=\\
\sin(\phi_j(t-\tau)-\phi_i(t))+\sin(\phi_j(t-\tau)+\phi_i(t)).
\end{array}
  \label{eq:f}
\end{align}
The equilibria $\phi^\pm$ in equation~\eqref{eq:fullphasemodel} are:
\begin{align}
  \begin{array}{l}
    \phi^+(n)=\dfrac{1}{2}\left(\arcsin\left(-\dfrac{1}{K}\right)+2n\pi\right)\\\\
    \phi^-(n)=\dfrac{1}{2}\left(\pi-\arcsin\left(-\dfrac{1}{K}\right)+2n\pi\right)\\
  \end{array},
\label{eq:equilibria}
\end{align}
$n\in\mathbb{Z},~K\geq 1$. For our analysis, we consider three main assumptions:
\begin{enumerate}[(a)]
\item\label{aa}The critical eigenvalue $\lambda$ of the linearization of~\eqref{eq:fullphasemodel} at equilibria crosses the imaginary axis with non vanishing velocity, i.e. $\re(\lambda'(\phi^\pm))\neq 0$.
\item\label{ab}The purely imaginary eigenvalue $\lambda=\i\omega$ is simple.
\item\label{ac}The linearization of~\eqref{eq:fullphasemodel} at equilibria, has no eigenvalues of the form $\i k \omega$, $k\in\Z-\{1,-1 \}$.
\end{enumerate}
The linerization of \eqref{eq:fullphasemodel} at equilibria is:
\begin{align}
\begin{array}{l}
\delta\ddot\phi_i+\mu\delta\dot\phi_i\\
-\dfrac{K\mu}{N-1}\displaystyle\sum_{\stackrel[j=1]{j\neq i}{}}^N\sum_{r=1}^\infty\bigg\{\dfrac{1}{r!}\bigg(\delta\phi_i\frac{\partial}{\partial\phi'_i}+\\
\delta\phi_{j\tau}\dfrac{\partial}{\partial\phi'_{j\tau}} \bigg)^rf(\phi_i,\phi_{j\tau}) \bigg\}_{\stackrel[\phi'_{j\tau}=\phi^\pm]{\phi'_i=\phi\pm}{}}=0.
\end{array}
\end{align}
where $\phi_{j\tau}:=\phi_j(t-\tau)$. Truncate the Taylor series up to the third-order term:
\note{DF: Computations:
  \begin{itemize}
  \item r=1
    \begin{align*}
      a_1= \left(\delta\phi_i\frac{\partial}{\partial\phi'_i}+\delta\phi_{j\tau}\frac{\partial}{\partial\phi'_{j\tau}}\right)f(\phi'_j,\phi'_{j\tau})\bigg|_{\stackrel[\phi'_{j\tau}=\phi^\pm]{\phi'_i=\phi\pm}{}}
    \end{align*}
    \begin{align*}
      \frac{\partial}{\partial\phi'_i}f(\phi'_i,\phi'_{j\tau})&=\cos(\phi'_{j\tau}-\phi'_i)(-1)+\cos(\phi'_{j\tau}+\phi'_i)\\
      \frac{\partial}{\partial\phi'_{j\tau}}f(\phi'_i,\phi'_{j\tau})&=\cos(\phi'_{j\tau}-\phi'_i)+\cos(\phi'_{j\tau}+\phi'_i)
    \end{align*}
    \begin{align*}
      a_1=\delta\phi_i(\cos 2\phi^\pm -1)+\delta\phi_{j\tau}(\cos 2\phi^\pm+1)
    \end{align*}
\end{itemize}
}
\note{
\begin{itemize}
\item r=2
  \begin{align*}
    a_2=\dfrac{1}{2!}\left(\delta\phi_i^2\frac{\partial^2}{\partial\phi_i^{'2}}+2\delta\phi_i\delta\phi_{j\tau}\frac{\partial^2}{\partial\phi'_i\partial\phi'_{j\tau}}+\delta\phi^2_{j\tau}\frac{\partial^2}{\partial_{j\tau}^{'2}} \right)f(\phi'_j,\phi'_{j\tau})\bigg|_{\stackrel[\phi'_{j\tau}=\phi^\pm]{\phi'_i=\phi\pm}{}}
  \end{align*}
  \begin{align*}
    \frac{\partial^2f}{\partial\phi_i^{'2}}&=\frac{\partial}{\partial\phi'_i}\left(\cos(\phi'_{j\tau}-\phi'_i)(-1)+\cos(\phi'_{j\tau}+\phi'_i) \right)\\
                                       &=\sin(\phi'_{j\tau}-\phi'_i)(-1)-\sin(\phi'_{j\tau}+\phi_i)\\
                                       &=-\sin 2\phi^\pm\\
\frac{\partial^2f}{\partial\phi_i^{'}\partial\phi'_{j\tau}}&=\frac{\partial}{\partial\phi'_i}\left(\cos(\phi'_{j\tau}-\phi'_i)+\cos(\phi'_{j\tau}+\phi'_i) \right)\\
                                                        &=-\sin(\phi'_{j\tau}-\phi'_i)(-1)-\sin(\phi'_{j\tau}+\phi'_i)\\
                                                        &=-\sin 2\phi^\pm\\
\frac{\partial^2f}{\partial\phi_{j\tau}^{'2}}&=\frac{\partial}{\partial\phi'_{j\tau}}\left(\cos(\phi'_{j\tau}-\phi'_i)+\cos(\phi'_{j\tau}+\phi'_i) \right)\\
                                          &=-\sin(\phi'_{j\tau}-\phi'_j)-\sin(\phi'_{j\tau}+\phi'_i)\\
                                          &=-\sin 2\phi^\pm
  \end{align*}
  \begin{align*}
    a_2=-\frac{1}{2}\sin 2\phi^\pm(\delta\phi_i+\delta\phi_{j\tau})^2
  \end{align*}
  \end{itemize}
}
\pagestyle{otherplain}
\note{
  \begin{itemize}
  \item r=3
    \begin{align*}
      a_3&=\frac{1}{3!}\bigg(\delta\phi_i^3\frac{\partial^3}{\partial\phi^{'3}_i}+3\delta\phi^2_i\delta\phi_{j\tau}\frac{\partial^3}{\partial\phi^{'2}_i\partial\phi'_{j\tau}}+3\delta\phi_i\delta\phi_{j\tau}^2\frac{\partial^3}{\partial\phi'_i\partial\phi^{'2}_{j\tau}}\\
         &+\delta\phi_{j\tau}^3\frac{\partial^3}{\partial\phi^{'3}_{j\tau}}\bigg)f(\phi'_j,\phi'_{j\tau})\bigg|_{\stackrel[\phi'_{j\tau}=\phi^\pm]{\phi'_i=\phi\pm}{}}
    \end{align*}
    \begin{align*}
      \frac{\partial^3f}{\partial\phi^{'3}_i}&=\frac{\partial}{\partial\phi'_i}\left(-\sin(\phi'_{j\tau}-\phi'_i)-\sin(\phi'_{j\tau}+\phi'_i) \right)\\
                                            &=-\cos(\phi'_{j\tau}-\phi'_i)(-1)-\cos(\phi'_{j\tau}+\phi'_i)\\
                                            &=1-\cos 2\phi^\pm\\
      \frac{\partial^3f}{\partial\phi^{'2}_i\partial\phi'_{j\tau}}&=\frac{\partial}{\partial\phi'_i}\left(\sin(\phi'_{j\tau}-\phi'_j)-\sin(\phi'_{j\tau}+\phi'_i) \right)\\
                                            &=\cos(\phi'_{j\tau}-\phi'_i)(-1)-\cos(\phi'_{j\tau}+\phi'_i)\\
                                            &=-(1+\cos 2\phi^\pm)\\
\frac{\partial^3f}{\partial\phi'_i\partial\phi^{'2}_{j\tau}}&=\frac{\partial}{\partial\phi'_i}\left(-\sin(\phi'_{j\tau}-\phi'_i)-\sin(\phi'_{j\tau}+\phi'_i) \right)\\
                                            &=-\cos(\phi'_{j\tau}-\phi'_i)(-1)-\cos(\phi'_{j\tau}+\phi'_i)\\
                                            &=1-\cos 2\phi^\pm\\
      \frac{\partial^3f}{\partial\phi^{'3}_{j\tau}}&=\frac{\partial}{\partial\phi'_{j\tau}}\left(-\sin(\phi'_{j\tau}-\phi'_i)-\sin(\phi'_{j\tau}+\phi'_i) \right)\\
                                           &=-\cos(\phi'_{j\tau}-\phi'_i)-\cos(\phi'_{j\tau}+\phi'_i)\\
                                           &=-(1+\cos 2\phi^\pm)
    \end{align*}
    \begin{align*}
      a_3=\frac{1}{6}\left[\left(\delta\phi_i-\delta\phi_{j\tau} \right)^3-\left(\delta\phi_i+\delta\phi_{j\tau} \right)^3\cos 2\phi^\pm\right]
    \end{align*}
  \end{itemize}
}
\begin{align}
\begin{array}{l}
  \ddot\phi_i+\mu\dot\phi_i=\\
\dfrac{K\mu}{N-1}\displaystyle\sum_{\stackrel[j=1]{j\neq i}{}}^N\bigg\{\left(\phi_{j\tau}-\phi_i \right)+\left(\phi_{j\tau}+\phi_i \right)\cos 2\phi^\pm\\
-\dfrac{1}{2}\left(\phi_{j\tau}+\phi_i \right)^2\sin 2\phi^\pm\\
-\dfrac{1}{6}\left[\left(\phi_{j\tau}-\phi_i \right)^3+\left(\phi_{j\tau}+\phi_i \right)^3\cos 2\phi^\pm \right]\bigg\}
\end{array}
\label{eq:linear_expansion_form}
\end{align}
here, for the sake of notation we changed $\delta\phi_i\to\phi_i$.

 Note that $\phi_i\in\kC([-\tau,0),\R)$, $i=1,\ldots,N$.
 Choosing $x^{(i)}_1=\phi_i$ and $x^{(i)}_2=\dot\phi_i$, the vector field form, is:
\begin{align}
\begin{array}{l}
\dot{x}_1^{(i)}=x_2^{(i)}\\
\dot{x}_2^{(i)}=-\mu x_2^{(i)}+\dfrac{K\mu}{N-1}\displaystyle\sum_{\stackrel[j=1]{j\neq i}{}}^N\bigg\{(-1+\cos 2\phi^\pm)x^{(i)}_1 \\
             ~~~+(1+\cos 2\phi^\pm)x^{(j)}_{1\tau} -\dfrac{1}{2}\left(x^{(j)}_{1\tau}+x^{(i)}_1 \right)^2\sin 2\phi^\pm\\
             ~~~-\dfrac{1}{6}\left[\left(x^{(j)}_{1\tau}-x^{(i)}_1 \right)^3+\left(x^{(j)}_{1\tau}+x^{(i)}_1 \right)^3\cos 2\phi^\pm \right]\bigg\}.
\end{array}
  \label{eq:vector_form}
\end{align}
We define $\kX:=\kC([-\tau,0],\R^{2N})$, the Banach space of continuous functions from $[-\tau,0]$ into $\R^{2N}$ equipped with the usual norm $$\|x\|=\stackrel[-\tau\leq\theta\leq0]{}{\text{sup}}|x(\theta)|,~~~x\in\kC([-\tau,0],\R^{2N}),$$
and $x=(x^{(1)},\ldots,x^{(N)})\in\kX$, where $x^{(i)}=(x_1^{(i)},x_2^{(i)})$.

Now, in order to build the decomposition of the infinite-dimensional space, we need to define de adjoint operator associated to the linear part of the linearization and a inner product, via a bilinear form.

Following~\cite{Hale1971, Gilsinn2009}, we can represent the dynamics in~\eqref{eq:vector_form} by the abstract differential equation:
\begin{align}
\frac{d}{dt} x_t(\phi)=A(\eta)x_t(\phi)+\kF(x_t(\phi),\eta),
  \label{eq:abstractODE}
\end{align}
which satisfies $(T(t)\phi)(\theta)=(x_t(\phi))(\theta)=x(t+\theta)$, where $T(t)$ is a semigroup of family of operators, $\theta\in[-\tau,0]$, and $\eta$ is a vector of parameters. The linear operator $A(\eta)\in\Mat(2N)$ is defined in equation $(2.16)$ in~\cite{FerruzzoCorrea2014a}, and
\begin{align}
  \left(\kF(x)\right)(\theta)=\left\{
    \begin{array}{ll}
      \frac{\partial x}{\partial\theta}(\theta)&,-\tau\leq\theta<0\\
      F(x(0),x(-\tau),\eta)&,\theta=0
    \end{array}
\right.,
\label{eq:general_F}
\end{align}
$F=(f^{(1)},\ldots,f^{(N)})^T$, where $f^{(i)}=(f^{(i)}_1,f^{(i)}_2)$, $f^{(i)}_1=0$, and $f^{(i)}_2=\dfrac{K\mu}{N-1}\displaystyle\sum_{\stackrel[j=1]{j\neq i}{}}^N\bigg\{-\dfrac{1}{2}\left(x^{(j)}_{1\tau}+x^{(i)}_1 \right)^2\sin 2\phi^\pm-\dfrac{1}{6}\left[\left(x^{(j)}_{1\tau}-x^{(i)}_1 \right)^3+\left(x^{(j)}_{1\tau}+x^{(i)}_1 \right)^3\cos 2\phi^\pm \right]\bigg\}$.

Associated to the linear part of~\eqref{eq:abstractODE}, we define the formal adjoint equation:
\begin{align}
  \dfrac{dy}{dt}(t,\eta)=A_0^T(\eta)y(t,\eta)+A_\tau^T(\eta)y(t+\tau,\eta),
\label{eq:adjoint}
\end{align}
matrices $A_0(\eta),~A_\tau(\eta)$  are defined in~\cite{FerruzzoCorrea2014a}. The strongly continuous semigroup $(T^*(t)\psi)(\theta)=(y_t(\psi))(\theta)=y(t+\theta)$, defines the infinitesimal generator:
\begin{align}
\begin{array}{l}
  (A^*(\eta)\psi)=\\
\left\{\begin{array}{ll}
    \frac{\partial\psi}{\partial\theta}(\theta)&,0<\theta\leq\tau\\
    A_0(\eta)^T\psi(0)+A_\tau(\eta)^T\psi(\tau)&,\theta=0
    \end{array}\right.,
\end{array}
\label{eq:adjoint_gen}
\end{align}
such that $\frac{d}{dt}T^*(t)\psi=A^*T^*(t)\psi$, $\psi\in\kX^*:=\kC([0,\tau],\R^{2N})$. The natural inner product, following~\cite{Hale1993}, has the form:
\begin{align*}
  \langle x,y \rangle=\bar{x}^T(0)y(0)+\int_{-\tau}^0\bar{x}(s+\tau)A_\tau(\eta)y(s)ds.
\end{align*}
For $\varphi\in\kX$ and $\psi\in\kX^*$ we have~\cite{Gilsinn2009}:
\begin{enumerate}
\item $\lambda$ is an eigenvalue of $A(\eta)$ if and only if $\bar\lambda$ is and eigenvalue of $A^*(\eta)$.
\item If $\varphi_1,\ldots,\varphi_d$ is a basis for the eigenspace of $A(\eta)$ and $\psi_1,\ldots,\psi_d$ is a basis for the eigenspace of $A^*(\eta)$, construct the matrices $\Phi=(\varphi_1,\ldots\varphi_d)$ and  $\Psi=(\psi_1,\ldots,\psi_d)$. Define the bilinear form:
  \begin{align}
    \langle \Psi,\Phi \rangle=I
\label{eq:bilinear_form}
  \end{align}
\end{enumerate}
\section{The Fixed Point space $\S_N$}
Due to the $\S_N$-symmetry of~\eqref{eq:fullphasemodel} the space where solutions $\phi_i$ lie can be decomposed into the Fixed-point subspace where symmetry-preserving solutions emerge  and a subspace with symmetry-breaking solutions, this was shown in~\cite{FerruzzoCorrea2014a}. We analyze stability of the periodic solutions near Hopf bifurcations in the Fixed point space, these bifurcations satisfy assumptions~\eqref{aa}-\eqref{ac} for $K>1$. In this subspace, equation~\eqref{eq:vector_form} has the form:
\begin{align}
  \begin{array}{l}
    \dot x_1=x_2\\
    \dot x_2=-\mu x_2+K\mu(-1+\cos 2\phi^\pm)x_1\\
            ~~~+K\mu(1+\cos 2\phi^\pm)x_{1\tau}-\dfrac{1}{2}(x_{1\tau}+x_1)^2\sin 2\phi^\pm\\
            ~~~-\dfrac{1}{6}\left[(x_{1\tau}-x_1)^3+(x_{1\tau}+x_1)^3\cos 2\phi^\pm \right],
  \end{array}
\end{align}
then matrices $A_0(\eta)$ and $A_\tau(\eta)$ in~\eqref{eq:adjoint_gen} become:
\begin{align}
A_0(\eta) = \left(
  \begin{array}{cc}
    0&1\\
K\mu(-1+\cos(2\phi^\pm))&-\mu
  \end{array}
\right),
\label{eq:A0}
\end{align}

\begin{align}
A_\tau(\eta)=\left(
  \begin{array}{cc}
    0&0\\
    K\mu(1+\cos(2\phi^\pm))&0
  \end{array}
\right),
  \label{eq:Atau}
\end{align}
and $F$ in~\eqref{eq:general_F} takes the form $F=(f_1~f_2)^T$, with $f_1=0$, and $f_2$:
\begin{align}
\begin{array}{l}
f_2(x_t,\eta)=-\dfrac{1}{2}(x_{1\tau}+x_1)^2\sin 2\phi^\pm\\
~~~~~~~~~~~~-\dfrac{1}{6}\left[(x_{1\tau}-x_1)^3+(x_{1\tau}+x_1)^3\cos 2\phi^\pm \right].
\end{array}
  \label{eq:F}
\end{align}
We need the complex eigenfunctions $As(\vartheta)=\i\omega s(\vartheta)$, $A^*n(\theta)=\i\omega n(\theta)$, associated to the critical eigenvalues $\lambda=\i\omega$, and $\bar \lambda=-\i\omega$ with $s(\vartheta) = s_1(\vartheta)+\i s_2(\vartheta)$ and $n(\theta)=n_1(\theta)+\i n_2(\theta)$. These eigenfunctions can be computed solving the boundary value problem $\frac{d}{d\vartheta}s_{1,2}=\mp\omega s_{2,1}(\vartheta)$, and $\frac{d}{d\theta}n_{1,2}=\pm\omega s_{2,1}(\vartheta)$, which, after substituting the operator $A(\eta)$, becomes:
\begin{align}
\begin{array}{rcl}
  A_0(\eta)s_1(0)+A_\tau(\eta)s_1(-\tau)& = &-\omega s_2(0)\\
  A_0(\eta)s_2(0)+A_\tau(\eta)s_2(-\tau)& = &\omega s_1(0)
\end{array}
\label{eq:A0_sol}
\end{align}
and
\begin{align}
\begin{array}{rcl}
  A_0^T(\eta)n_1(0)+A_\tau^T(\eta)n_1(-\tau)& = &\omega n_2(0)\\
  A_0^T(\eta)n_2(0)+A_\tau^T(\eta)n_2(-\tau)& = &-\omega n_1(0)
\end{array},
\label{eq:Atau_sol}
\end{align}
with general solutions:
\begin{align}
\begin{array}{rcl}
s_1(\vartheta) &=& \cos(\omega\vartheta)c_1-\sin(\omega\vartheta)c_2\\
s_2(\vartheta) &=& \sin(\omega\vartheta)c_1+\cos(\omega\vartheta)c_2\\
n_1(\theta)&=& \cos(\omega\theta)d_1-\sin(\omega\theta)d_2\\
n_2(\theta) &=& \sin(\omega\theta)d_1+\cos(\omega\theta)d_2
\end{array}.
\label{eq:gen_sol}
\end{align}
The coefficients $c_1=[c_{11}~c_{12}]^T,~c_2=[c_{21}~c_{22}]^T,~d_1=[d_{11}~d_{12}]^T,~d_2=[d_{21}~d_{22}]^T$ can be obtained by considering the boundary conditions,
\begin{align}
\begin{array}{l}
  \left(\begin{array}{c}
    A_0(\eta)+\cos(\omega\tau)A_\tau(\eta)\\\omega I+\sin(\omega\tau)A_\tau(\eta)
  \end{array}\right)^T\left(
  \begin{array}{c}
    c_1\\c_2
  \end{array}
\right)=0\\\\
 \left( \begin{array}{c}
    A_0^T(\eta)+\cos(\omega\tau)A_\tau^T(\eta)\\
-\omega I-\sin(\omega\tau)A_\tau^T(\eta)
  \end{array}\right)^T\left(
  \begin{array}{c}
    d_1\\d_2
  \end{array}
\right)=0
\end{array},
\end{align}
the ``orthonormality'' condition $\langle s,n\rangle=I$, and setting $c_{11}=1$ and $c_{21}=0$, see ~\cite{Kalmar-Nagy2001, Hale1977} for more details.

It is also possible to decompose the solution $x_t(\vartheta)$ to equation~\eqref{eq:abstractODE} into $x_t(\vartheta)=y_1(t)s_1(\vartheta)+y_2(t)s_2(\vartheta)+\w(t)(\vartheta)$, where $y_1$ and $y_2$ lie in the center subspace, such that $y_{1,2}(t)=\langle n_{1,2}(0),x_t(0) \rangle$, and $\w$ in the infinite-dimensional component subspace, thus, we have
\begin{align}
  \begin{array}{rcl}
    \dot y_1 &=& \omega y_2 + n_1^T(0)F\\
    \dot y_2 &=& -\omega y_1 + n_2^T(0)F\\
\end{array}
\label{eq:y1y2}
\end{align}
\begin{align}
     \dot \w  = A(\eta)+\mathcal{F}(x_t,\eta)-n_1^T(0)Fs_1 - n_2^T(0)Fs_2,
\label{eq:trans}
\end{align}
where
\begin{align}
\begin{array}{l}
  \mathcal{F}=\\
\left\{\begin{array}{l}
      0,~~~\vartheta\in[-\tau,0)\\
      F(y_1(t)s_1(0)+y_2(t)s_2(0)+\w(t)(0)),~\vartheta=0
    \end{array}\right..
\end{array}
\label{eq:FF-F}
\end{align}
\subsection{The center manifold}
Following~\cite{Hassard1981, Kalmar-Nagy2001, Zhao2009}, we know that $\w$ can be approximated by the second-order expansion:
\begin{align}
 \w(y_1,y_2)(\vartheta)=\dfrac{1}{2}(h_1(\vartheta)y_1^2+2h_2(\vartheta)y_1y_2+h_3(\vartheta)y_2^2),
\label{eq:center_manifold}
\end{align}
thus, by differentiating and substituting equation~\eqref{eq:trans} keeping up to second order terms, we obtain:
\begin{align}
  \dot\w = -\omega h_2y_1^2+\omega(h_1-h_3)y_1y_2 + \omega h_2 y_2^2  + O(y^3),
\label{eq:dot_w}
\end{align}
and from equation~\eqref{eq:trans}, 
\begin{align}
\begin{array}{l}
  \dfrac{d\w}{dt}=\\
A(\eta)\w+\mathcal{F}(\w+y_1s_1+y_2s_2)-(d_{12}s_1+d_{22}s_2)f_2.
\end{array}
\label{eq:dwdt}
\end{align}
From the definition of $A(\eta)$, equivalent to~\eqref{eq:adjoint_gen}, we see that
\begin{align}
\begin{array}{l}
  A(\eta)\w=\\
\left\{\begin{array}{ll}
      \frac{1}{2}(\dot h_1 y_1^2+2\dot h_2 y_1y_2 + \dot h_3y_2^2),&\vartheta\in[-\tau,0)\\
      A_0(\eta)\w(0)+A_\tau(\eta)\w(-\tau),&\vartheta=0
    \end{array}\right .,
\end{array}
\label{eq:Aw}
\end{align}
then, from equation~\eqref{eq:center_manifold},~\eqref{eq:dot_w},~\eqref{eq:dwdt}, and~\eqref{eq:Aw}, we can obtain the unknown coefficients $h_1,~h_2$, and $h_3$, solving:
\begin{align}
  \begin{array}{rcl}
    \dot h_1 &=& 2( -\omega h_2 +f_2^{20}(d_{12}s_1(\vartheta)+d_{22}s_2(\vartheta))),\\
    \dot h_2 &=& \omega(h_1-h_3)+ f_2^{11}(d_{12}s_1(\vartheta)+d_{22}s_2(\vartheta)),\\
    \dot h_3 &=& 2(\omega h_2 + f_2^{02}(d_{12}s_1(\vartheta)+d_{22}s_2(\vartheta))),\\
  \end{array}
\label{eq:h1h2h3a}
\end{align}
\begin{align}
  \begin{array}{l}
A_0(\eta)h_1(0)+ A_\tau(\eta) h_1(-\tau)=\\
~~~~~~~~~2(-\omega h_2(0)+ f_2^{02}(d_{12}s_1(0)+d_{22}s_2(0))),\\\\
A_0(\eta)h_2(0)+A_0 h_2(-\tau)=\\
~~~\omega (h_1(0)-h_3(0)) + f_2^{11}(d_{12}s_1(0)+d_{22}s_2(0))),\\\\
A_0(\eta)h_3(0)+ A_\tau(\eta) h_3(-\tau)=\\
~~~~~~~~~2(\omega h_2(0)+ f_2^{02}(d_{12}s_1(0)+d_{22}s_2(0))),
  \end{array}
\label{eq:h1h2h3}
\end{align}
where $f^{20} = \dfrac{1}{2}\dfrac{\partial^2 f}{\partial y_1^2}\bigg|_0$, $f^{11} = \dfrac{\partial^2 f}{\partial y_1\partial y_2}\bigg|_0$, and $f^{02} = \dfrac{1}{2}\dfrac{\partial^2 f}{\partial y_2^2}\bigg|_0$.

Equation~\eqref{eq:h1h2h3a} is written as the inhomogeneus differential equation:
\begin{align}
  \frac{dh}{d\vartheta}= Ch+p\cos(\omega\vartheta)+q\sin(\omega\vartheta)
\label{eq:inhomogeneus_h}
\end{align}
where
\begin{align*}
\begin{array}{l}
  h:= \left(
    \begin{array}{c}
      h_1\\h_2\\h_3
    \end{array}
\right),~~~C:= \omega\left(
  \begin{array}{ccc}
    0&-2 I &0\\
    I&0&-I\\
    0&2 I&0
  \end{array}
\right)_{6\times 6}\\\\
p:=\left(
  \begin{array}{c}
    f_2^{02}p_0\\f_2^{11}p_0\\f_2^{02}p_0
  \end{array}
\right),~~~q:=\left(
  \begin{array}{c}
    f_2^{02}q_0\\f_2^{11}q_0\\f_2^{02}q_0
  \end{array}
\right),\\\\
p_0:=\left(
  \begin{array}{c}
    d_{12}\\c_{22}d_{22}
  \end{array}
\right),~~~q_0:=\left(
  \begin{array}{c}
    d_{22}\\-c_{22}d_{12}
  \end{array}
\right),
\end{array}
\end{align*}
with general solution:
\begin{align}
h(\vartheta)=\e^{C\vartheta}K + M\cos(\omega\vartheta) + N\sin(\omega\vartheta).
  \label{eq:gen_sol}
\end{align}
After substituting the general solution  into~\eqref{eq:inhomogeneus_h} we solve for $M$ and $N$, and then from the boundary value problem we solving for $K$,
 \begin{align}
\left(
  \begin{array}{cc}
    C&-\omega I\\
   \omega I& C
  \end{array}
\right)\left(
  \begin{array}{c}
    M\\N
  \end{array}
\right)=-\left(
  \begin{array}{c}
    p\\q
  \end{array}
\right)
   \label{eq:MN}
 \end{align}
 \begin{align}
   Ph(0)+Qh(-\tau)=p-r,
\label{eq:K_sol}
 \end{align}
where
\begin{align}
\begin{array}{l}
  P:=\left(
    \begin{array}{ccc}
      A_0&0&0\\
      0&A_0&0\\
      0&0&A_0
    \end{array}
\right)-C,\\\\
Q:=\left(
    \begin{array}{ccc}
      A_\tau&0&0\\
      0&A_\tau&0\\
      0&0&A_\tau
    \end{array}
\right),
\end{array}
\end{align}
and $r:=\left(
  \begin{array}{cccccc}
    0&f_2^{20}&0&f_2^{11}&0&f_2^{02}
  \end{array}
\right)^T$.

The expressions for $\w_1(0)$ and $\w_1(-\tau)$, necessary in~\eqref{eq:FF-F}, are:
\begin{align}
\begin{array}{l}
  \w_1(0)=\dfrac{1}{2}\bigg( (M_1+K_1)y_1^2 + 2(M_3+K_3)y_1y_2\\
~~~~~~~~~+(M_5+K_5)y_2^2\bigg),\\\\
 \w_1(-\tau)=\dfrac{1}{2}\bigg((\e^{-C\tau}K|_1+M_1\cos(\omega\tau)\\
~~~~~~~~~-N_1\sin(\omega\tau))y_1^2+2(\e^{-C\tau}K|_3\\
~~~~~~~~~+M_3\cos(\omega\tau)-N_3\sin(\omega\tau))y_1y_2\\
~~~~~~~~~+(\e^{-C\tau}K|_5+M_5\cos(\omega\tau)\\
~~~~~~~~~-N_5\sin(\omega\tau))y_2^2\bigg),
\end{array}
\label{eq:W0Wtau}
\end{align}
note that we only need $\w_1(\vartheta)$ since the nonlinear function in~\eqref{eq:F} only depends on $x_1$; then by substituting~\eqref{eq:W0Wtau} into~\eqref{eq:y1y2}, we obtain: 
\begin{align}
\begin{array}{rr}
\dot y_1 =& \omega y_2 + g_1(y_1,y_2;\eta)\\
\dot y_2 =& -\omega y_1 + g_2(y_1,y_2;\eta)
\end{array},
\end{align}
or
\begin{align}
  \begin{array}{l}
    \dot y_1=\omega y_2 + a_{20}y_1^2 +a_{11}y_1y_2 + a_{02}y_2^2 + a_{30}y_1^3\\
~~~~~~+a_{21}y_1²y_2+a_{12}y_1y_2^2+ a_{03}y_2^3,\\
    \dot y_2=-\omega y_1 + b_{20}y_1^2 +b_{11}y_1y_2 + b_{02}y_2^2 + b_{30}y_1^3\\
~~~~~~+b_{21}y_1²y_2+b_{12}y_1y_2^2+ b_{03}y_2^3.
  \end{array}
\label{eq:doty1y2}
\end{align}
In \cite{Guckenheimer1983} is computed the coefficient $a$, which determines stability of the normal form~\eqref{eq:doty1y2},  
\begin{align}
\begin{array}{l}
  a=\dfrac{1}{16}\left[g^{03}_2+g^{21}_2+g^{12}_1+g^{30}_1 \right]\\
~~~~~~+\dfrac{1}{16\omega}\bigg[g^{11}_2\left(g^{02}_2+g^{20}_2 \right)-g^{11}_1\left(g^{02}_1+g^{20}_1 \right)\\
~~~~~~- g^{02}_2g^{02}_1+g^{20}_2g^{20}_1\bigg],
\end{array}
\label{eq:a}
\end{align}
where $g^{ij}_r=\dfrac{\partial^{i+j}}{\partial^iy_1\partial^jy_2}g_r(0,0)$. Periodic orbits near Hopf bifurcation at the critical eigenvalue $\lambda=\i\omega$, will be stable if $a<0$ and unstable if $a>0$.
\section{Numerical Results}
\label{sec:numres}
\begin{figure}[!htb]
\centering
  \includegraphics[width=0.5\textwidth]{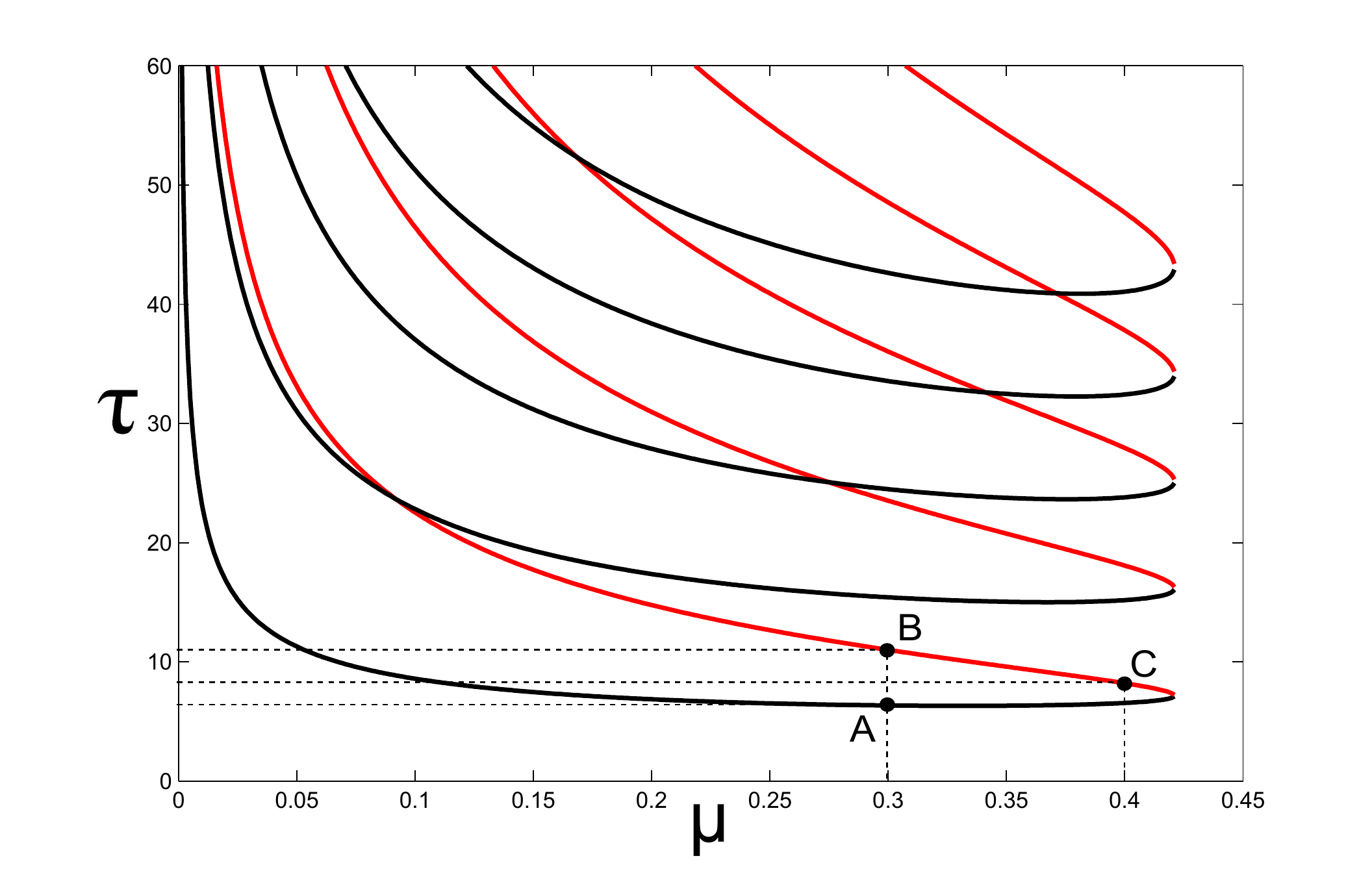}
  \caption{Symmetry-preserving bifurcations curves in $\Fix(\S_N)$ for $K=1.05$. In black, bifurcations crossing from the left to the right, in red, bifurcations crossing from the right to the left.}
  \label{fig:symm_bif_Fix}
\end{figure}

We reproduced some of the computations for the Hopf bifurcations in the Fixed point space for the case $K>1$ presented in~\cite{FerruzzoCorrea2014a}, we will compute stability for these bifurcation curves using results obtained in the previous section. In figure~\ref{fig:symm_bif_Fix} (part of figure 10, in~\cite{FerruzzoCorrea2014a}) are shown the symmetry-preserving bifurcations curves in the parameter space $(\mu,\tau)$ for $K=1.05$, for both cases: bifurcations crossing from the left to the right in black color, and crossing from the right to the left in red color; we  also choose three testing point for numerical simulation $A=(0.3,6.34)$, $B=(0.3,11)$, and $C=(0.4,8.204)$.

\begin{figure}[!htb]
\centering
  \includegraphics[width=0.5\textwidth]{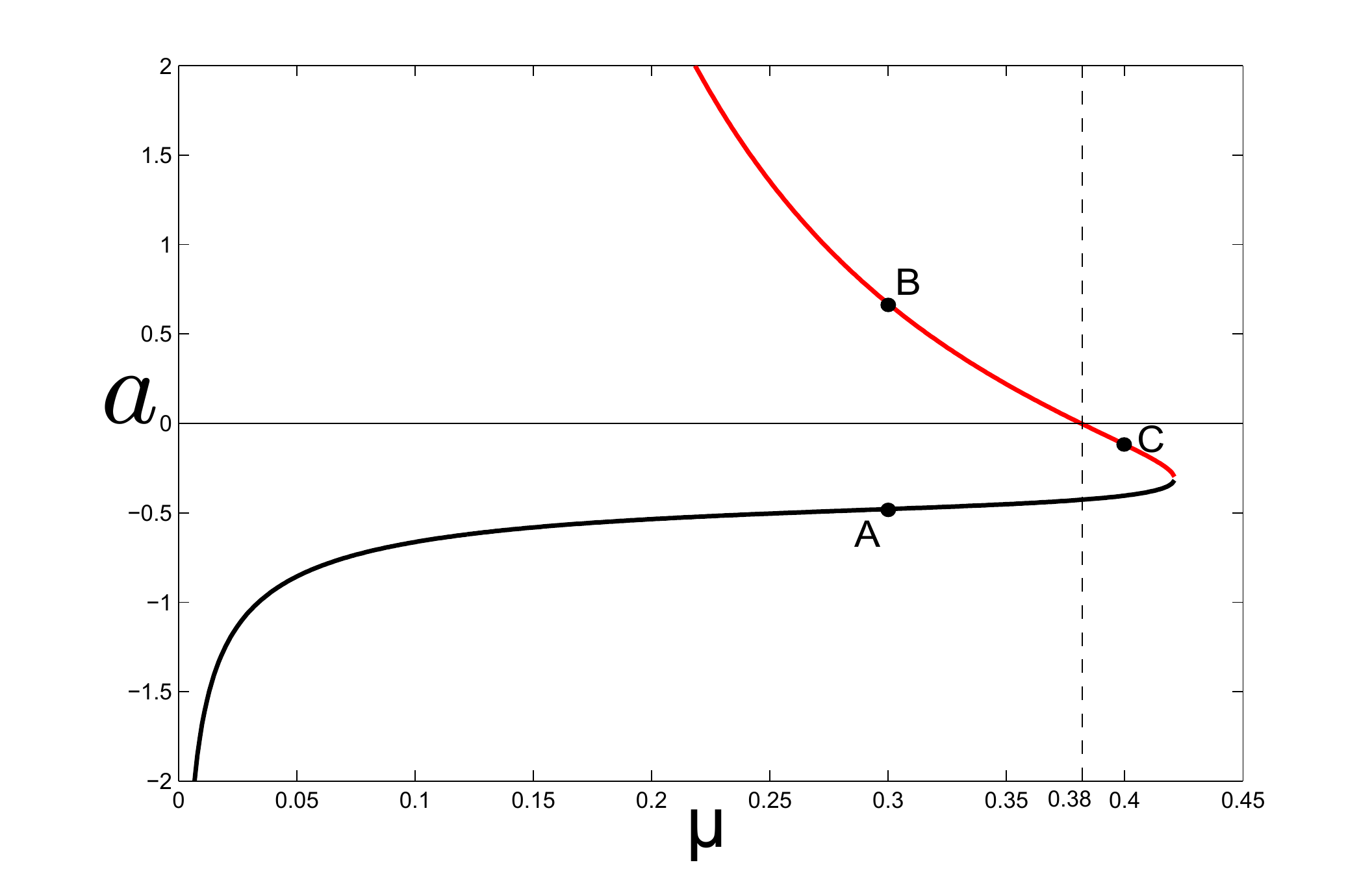}
\caption{Coefficient $a$ computed using eq.\eqref{eq:a}, determining stability of Hopf bifurcations in $\Fix(\S_N)$ for $K=1.05$, see figure~\ref{fig:symm_bif_Fix}. }
\label{fig:a}
\end{figure}

In figure~\ref{fig:a} is shown the coefficient $a$ computed using equation~\eqref{eq:a}, for $K=1.05$ in the parameter space $(\mu,\tau)$ related to the Hopf bifurcations shown in figure~\ref{fig:symm_bif_Fix}, the  black curve correspond to stability of Hopf bifurcations crossing from the left to the right (black curves in figure~\ref{fig:symm_bif_Fix}), as we can see, periodic solutions emerging at these Hopf bifurcations are all stable ($a<0$); the red curve correspond to stability of periodic orbits near Hopf bifurcations crossing back from the right to the left (red curves in figure~\ref{fig:symm_bif_Fix}), they are all unstable for $\mu<\mu^*(K)\approx 0.382$, and stable for $\mu>\mu^*$.

\begin{figure}[!htb]
\centering
  \includegraphics[width=0.5\textwidth]{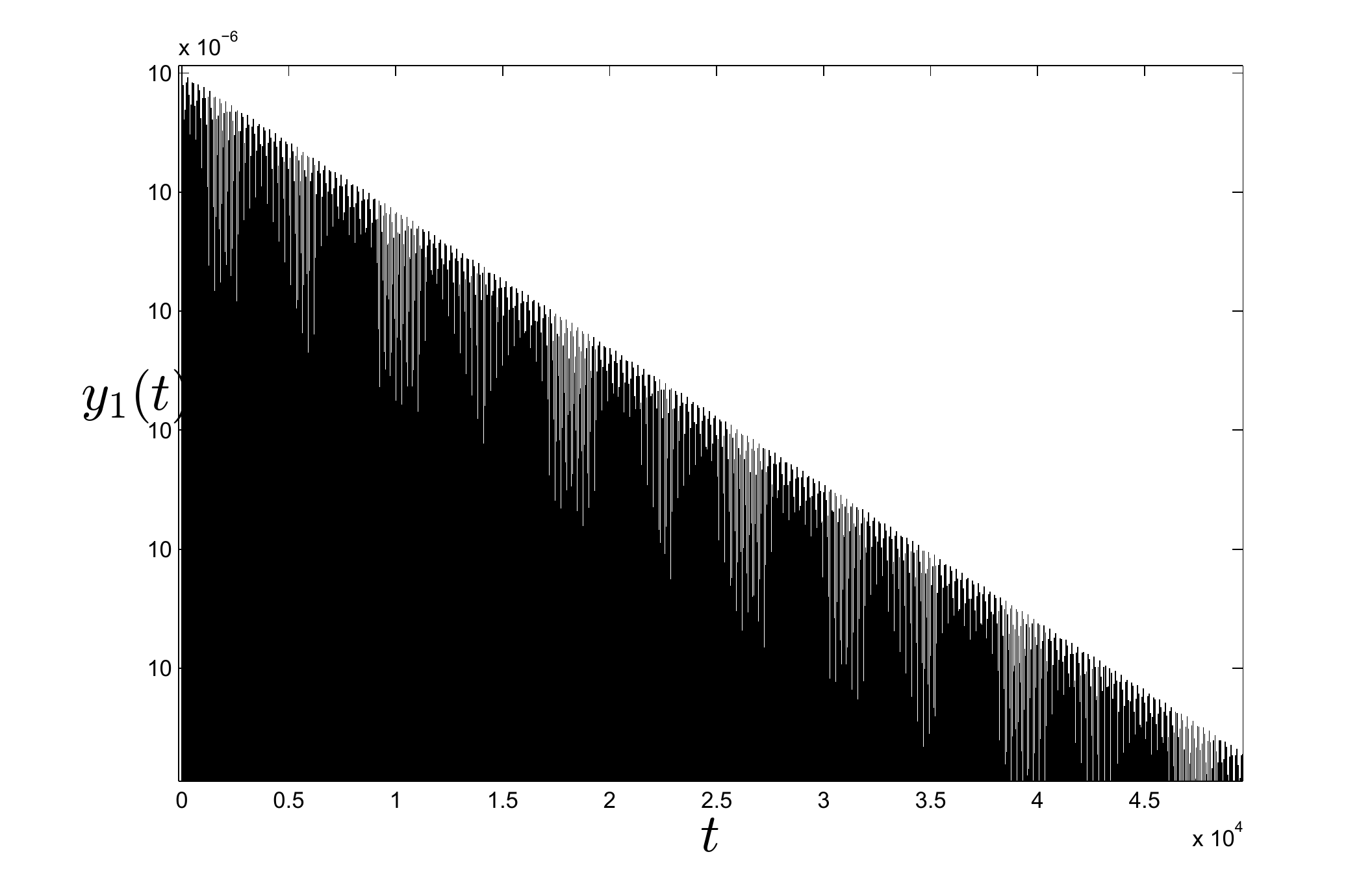}
  \caption{$y_1(t)$ at point $A$, for $\mu=0.3$ and $\tau=6.34$, c.i. $y_1(0)=y_2(0)=10^{-5}$. }
  \label{fig:y1_a}
\end{figure}

\begin{figure}[!htb]
\centering
  \includegraphics[width=0.5\textwidth]{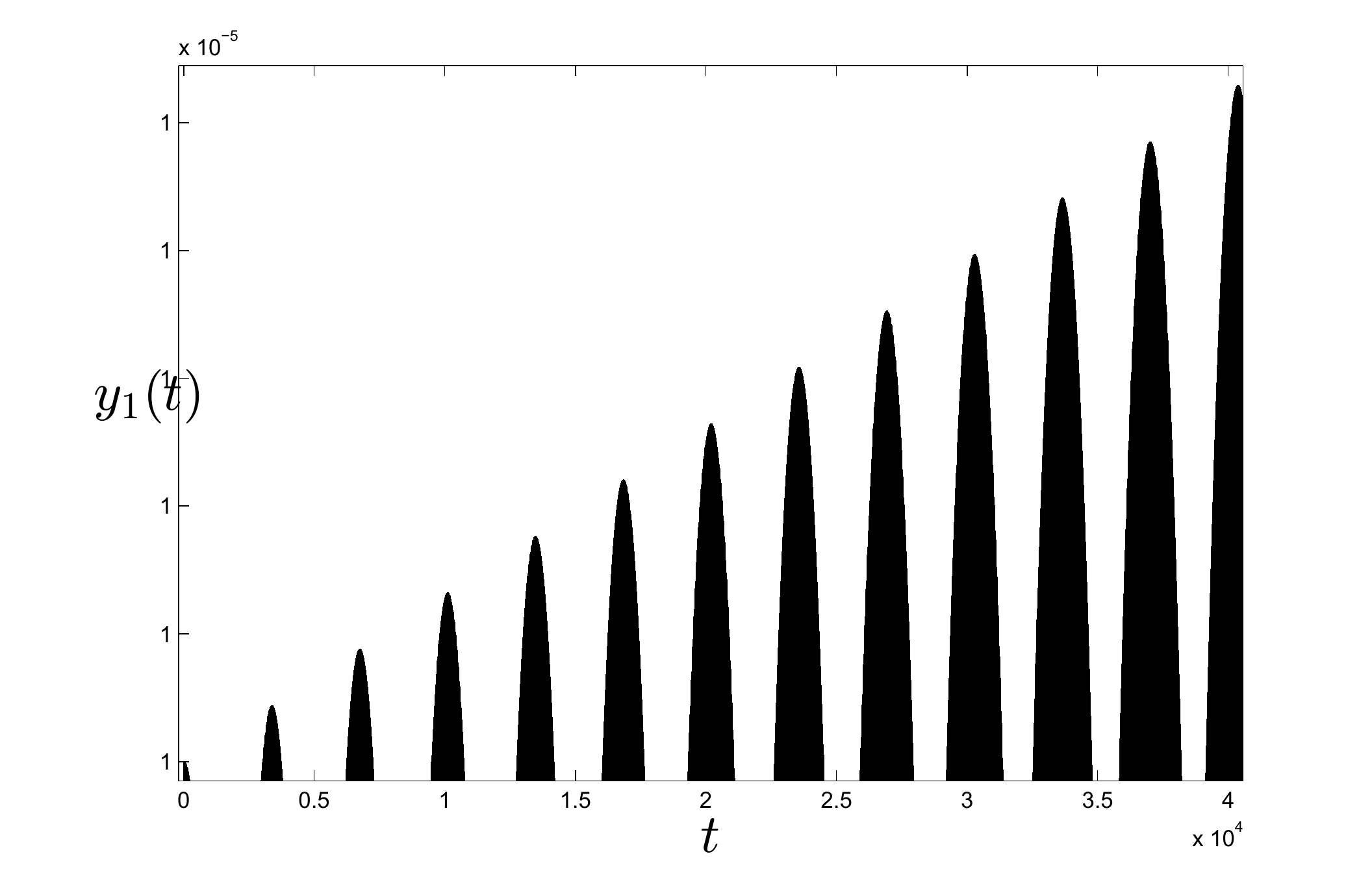}
\caption{$y_1(t)$ at point $B$, for $\mu=0.3$ and $\tau=11$, c.i. $y_1(0)=y_2(0)=10^{-5}$. }
\label{fig:y1_b}
\end{figure}

\begin{figure}[!htb]
\centering
  \includegraphics[width=0.5\textwidth]{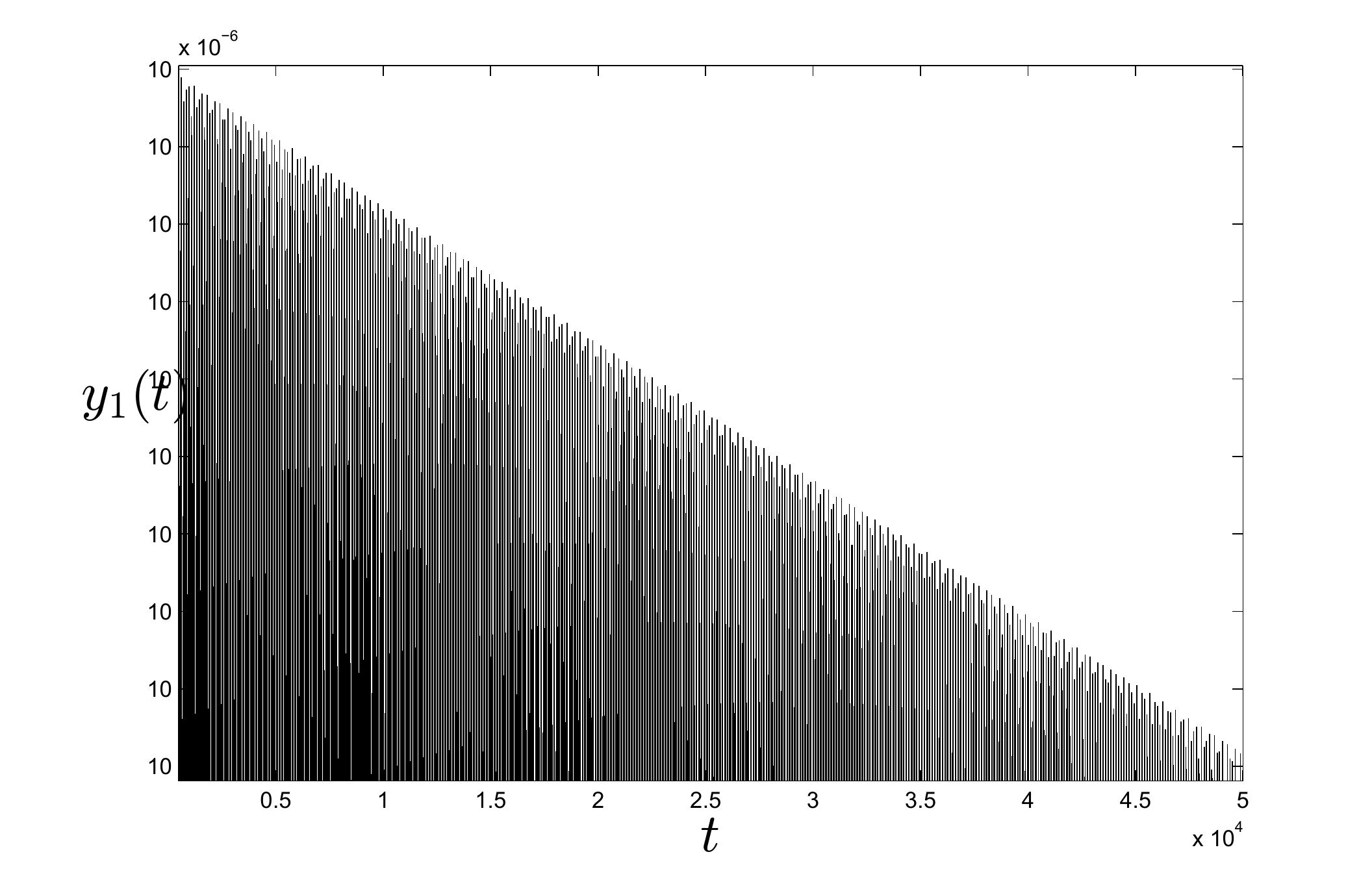}
\caption{$y_1(t)$ at point $C$, for $\mu=0.4$ and $\tau=8.204$, c.i. $y_1(0)=y_2(0)=10^{-5}$. }
\label{fig:y1_c}
\end{figure}

In order to confirm our results, numerical simulations computing $y_1$ and $y_2$ in equation~\eqref{eq:doty1y2} were run for point $A$, $B$, and $C$,  using \textit{ODE45} Matlab rutine, with time spam $5\times 10^4$, and maximum step size $0.05$. Periodic solution $y_1(t)$, for point $A$, $\mu=0.3$ and $\tau=6.34$, corresponding to Hopf bifurcation crossing from the left to the right in figure, which is stable ($a<0$), is  shown in figure~\ref{fig:y1_a}; periodic solution near Hopf bifurcations crossing from the right to the left at point $B$, $\mu=0.3$, $\tau=11$, which is unstable ($a>0$), is shown in figure~\ref{fig:y1_b}, and periodic solution at point $C$, $\mu=0.4$, $\tau=8.204$, which is stable ($a<0$), is shown in figure~\ref{fig:y1_c}; all initial conditions were set $y_1(0)=y_2(0)=10^{-5}$.
\section{Conclusions}
\label{sec:Conclusions}
The reduction of the inifinite-dimensional space onto the center manifold in normal form, was applied to the Fixed point space for the Full-phase model in order to analyse the stability of simple Hopf bifurcations, in both cases, for bifurcations crossing from the left to the right and for bifurcations crossing in the other way round. We found that near all bifurcations that cross from the stable region into the unstable region can emerge  periodic orbits which are stable ($a<0$), and, on the other hand, we observed that around all bifurcations coming back from the right to the left, unstable ($a>0$) periodic orbits can emerge for $\mu<\mu^*(K)$, and stable periodic orbits for $\mu>\mu^*(K)$.

Although, we computed the coefficient $a$ for a specific value of $K$, the procedure shown in this work is valid for all the parameter space where simple Hopf bifurcations appear. 

Finally, it is important to spotlight some points for further research: First, what is the nature of the solutions at the special point $\mu=\mu^*(K)$, at which the coefficient $a$ changes sign. Second, analyze stability of the degenerate Hopf bifurcations at the  Fixed point space for $K=1$, which are codimension 2, pure imaginary eigenvalue and zero eigenvalue, and third, the stability of the symmetry-breaking  degenerate Hopf  bifurcations which have multiplicity $N-1$.
\section*{Acknowledgements}
We would like to thank the Escola Polit\'ecnica da Universidade de S\~ao Paulo and FAPESP for their support.
\bibliographystyle{plain}
\bibliography{/home/diego/Dropbox/bibliography/phD_bibliography.bib}
\end{document}